\documentclass[11pt,twoside]{amsart} 
\usepackage{amssymb,pb-diagram,graphicx}
\advance \topmargin by -\headheight
\advance \topmargin by -\headsep
\evensidemargin \oddsidemargin
\newtheorem{theorem}{Theorem}[section] 
\newtheorem{corollary}[theorem]{Corollary}
\newtheorem{lemma}[theorem]{Lemma}
\newtheorem{proposition}[theorem]{Proposition}
\theoremstyle{definition} 
\newtheorem{definition}[theorem]{Definition}

\newtheorem{remark}[theorem]{Remark}
\newcommand{\Gm}{\mathbb G_{\mathrm{m}}}
\newcommand{\mO}{\mathcal O}
\title{Gersten resolutions for logarithmic forms} 
\author{M.Rovinsky} 
\address{AG Laboratory, HSE University, 6 Usacheva str., Moscow, Russia, 119048} 
\email{marat@mccme.ru}
\date{}
\begin{document} 
\begin{abstract} Let $X$ be a smooth variety over a field $k$ of characteristic 0. 

We show that (two versions of) the Zariski sheaves of regular $k$-linear logarithmic differential 
forms on $X$ (in the most na\"{\i}ve sense, i.e. sums of products of logarithmic differentials) 
have Gersten property (and therefore coincide). It turns out that the values of these sheaves 
are free abelian groups. As a consequence, for any integer $n>0$, these sheaves modulo 
$n$ also have Gersten property. The proof is an elementary application of the Voevodsky's 
$(-1)$-machinery and consists of assembling several well-known results. \end{abstract} 

\maketitle

\section{Notation and the statements} 

\subsection{} 
Let $k$ be a field of characteristic 0. For each $k$-scheme $X$ and each integer 
$s\ge 0$, denote by $X^{(s)}$ the set of points of $X$ of codimension $s$; for each 
$x\in X^{(s)}$, $k(x)$ denotes the residue field of $x$, and $i_x$ is the inclusion 
$x\hookrightarrow X$. In particular, for each abelian group $A$, $i_{x*}A$ is the 
constant sheaf on the closure of $x$ with the fibre $A$ extended to $X$ by zero. 

For each Zariski sheaf $\mathcal F$ (at least) on $X$, and each subscheme 
$Y\subseteq X$, set $\mathcal F(Y):=\mathop{\underrightarrow{\lim}}\limits_U\mathcal F(U)$, 
where $U$ runs over the open subsets of the closure of $Y$ containing $Y$. 

\subsection{} 
For each integer $q\ge 0$, denote by $\Omega^q_{-|k}$ the Zariski sheaf $U\mapsto\Omega^q_{U|k}(U)$ 
of regular $k$-linear differential $q$-forms on the category of smooth $k$-varieties. 

\subsubsection{$\Omega^{q,\log}_{-|k}$} 
For each smooth connected variety $U$ over $k$, define the subgroup 
$\Omega^{q,\log}_{-|k}(U)\subseteq\Omega^q_{-|k}(U)$ as the group of those regular $k$-linear 
$q$-forms on $U$ that are (finite) sums of exterior products of 
logarithmic differentials of non-zero rational functions on $U$: 
\[\Omega^{q,\log}_{-|k}(U):=\left\{\eta\in\Omega^q_{-|k}(U)\ |\ \eta=\sum_i
\frac{\mathrm{d}f_{i1}}{f_{i1}}\wedge\cdots\wedge\frac{\mathrm{d}f_{iq}}{f_{iq}}\mbox{ for some }
f_{ij}\in k(U)^{\times}\right\}\mbox{ if }q>0,\] and $\Omega^{0,\log}_{-|k}(U):=\mathbb Z$. 

\subsubsection{$\Omega^{q,\log}_{M,-|k}$} 
Set $\Omega^{0,\log}_{M,-|k}(U):=\mathbb Z$, and for $q>0$ let 
$\Omega^{q,\log}_{M,-|k}\subseteq\Omega^q_{-|k}$ be Zariski sheafification of the presheaf 
$U\mapsto\left\{\sum_i\frac{\mathrm{d}f_{i1}}{f_{i1}}\wedge\cdots\wedge
\frac{\mathrm{d}f_{iq}}{f_{iq}}\in\Omega^q_{-|k}(U)\ |\ f_{ij}\in\mathcal O(U)^{\times}\right\}$, 
so \[\Omega^{q,\log}_{M,-|k}(U)=\left\{\eta\in\Omega^{q,\log}_{-|k}(U)\ |\ \eta=\sum_i
\frac{\mathrm{d}f_{i1}^{(x)}}{f_{i1}^{(x)}}\wedge\cdots\wedge\frac{\mathrm{d}f_{iq}^{(x)}}{f_{iq}^{(x)}}
\begin{tabular}{l} for all $x\in U$ and\\ some $f_{ij}^{(x)}\in\mathcal O_{U,x}^{\times}$
\end{tabular}\right\}.\] 
In other words, $\Omega^{q,\log}_{M,-|k}$ is the image in $\Omega^q_{-|k}$ of the 
Zariski sheafification $\mathcal K^M_q$ of the presheaf $U\mapsto K^M_q(\mO(U))$ 
of $q$-th Milnor $K$-groups under the map $\mathrm{d}\log:\{f_1,\dots,f_q\}
\mapsto\frac{\mathrm{d}f_1}{f_1}\wedge\cdots\wedge\frac{\mathrm{d}f_q}{f_q}$. 

Then $\Omega^{1,\log}_{-|k}(U)=\Omega^{1,\log}_{M,-|k}(U)$ is canonically isomorphic to the 
quotient of the group $\mathcal O(U)^{\times}$ by the multiplicative group of the algebraic 
closure of $k$ in $k(U)$; for each point $x\in U$ and $j>0$, one has $\Omega^{j,\log}_{-|k}(x)
=\Omega^{j,\log}_{M,-|k}(x)=\left\{\sum_i\frac{\mathrm{d}f_{i1}}{f_{i1}}\wedge\cdots\wedge
\frac{\mathrm{d}f_{ij}}{f_{ij}}\ |\ f_{is}\in k(x)^{\times}\right\}\subset\Omega^j_{k(x)|k}$. 

\subsection{} 
We are going to deduce from some classical results the following one (probably known to the experts): 
\begin{theorem} \label{log} Let $X$ be a smooth $k$-variety, and 
$q\ge 0$ be an integer. Then the restriction of $\Omega^{q,\log}_{-|k}$ to $X$ with Zariski topology is 
a sheaf admitting a canonical flabby resolution 
\[\coprod_{x\in X^{(0)}}\!i_{x*}\Omega^{q,\log}_{-|k}(x)\to
\!\coprod_{x\in X^{(1)}}\!i_{x*}\Omega^{q-1,\log}_{-|k}(x)\to\cdots\to\!
\coprod_{x\in X^{(q-1)}}\!i_{x*}\Omega^{1,\log}_{-|k}(x)\to\!\coprod_{x\in X^{(q)}}
\!i_{x*}\mathbb Z\to 0.\] \end{theorem} 

\begin{corollary}[Bloch's formula] Let $X$ be a smooth variety over a field $k$ of characteristic 
$0$, and $q\ge 0$ an integer. Then $H^q(X,\Omega^{q,\log}_{X|k})$ is canonically isomorphic to 
the group $CH^q(X)$ of codimension $q$ cycles on $X$ modulo rational equivalence. \end{corollary} 

\begin{corollary} \label{free} For any pair $U\subseteq V$ of intersections of open 
subsets in a smooth variety over a field $k$ of characteristic $0$ and an integer 
$q\ge 0$, the restriction $\Omega^{q,\log}_{-|k}(V)\to\Omega^{q,\log}_{-|k}(U)$ is 
a spit inclusion of free abelian groups. \end{corollary}

\begin{corollary} \label{log-mod-n} Let $X$ be a smooth $k$-variety, $q\ge 0$ and $n>0$ be integers. 
Then the Zariski sheaf $\Omega^{q,\log}_{-|k}/n$ admits a canonical flabby resolution\\ 
\scalebox{0.93}{$\coprod\limits_{x\in X^{(0)}}\!i_{x*}\Omega^{q,\log}_{-|k}(x)/n\to
\!\coprod\limits_{x\in X^{(1)}}\!i_{x*}\Omega^{q-1,\log}_{-|k}(x)/n\to\cdots\to\!
\coprod\limits_{x\in X^{(q-1)}}\!i_{x*}\Omega^{1,\log}_{-|k}(x)/n\to\!\coprod\limits_{x\in X^{(q)}}
\!i_{x*}\mathbb Z/n\to 0$.} \end{corollary} 

One more corollary is the following `moving lemma': 
\begin{theorem} \label{moving_lemma} The sheaves $\Omega^{q,\log}_{-|k}$ 
and $\Omega^{q,\log}_{M,-|k}$ coincide. \end{theorem} 

The following result is a consequence of the Bloch--Kato conjecture proved in \cite{V_Bloch-Kato}. 
\begin{corollary} \label{alg-closed-mod-n} If $k$ algebraically closed then, for 
any integer $q\ge 0$, the sheaf $\Omega^{q,\log}_{-|k}$ is the quotient of the 
sheaf $\mathcal K^M_q$ by the subsheaf of divisible elements. \end{corollary} 
In particular, the norm residue isomorphism 
$\mathcal K^M_{q,X}/n\xrightarrow{\sim}\mathcal H^{q,\mathrm{\acute{e}t}}_X(\mu_n^{\otimes q})$ 
factors through isomorphisms $\mathcal K^M_{q,X}/n\xrightarrow[\raisebox{0.4ex}{$\sim$}]{\mathrm{d}\log_n}
\Omega^{q,\log}_{X|k}/n\xrightarrow{\sim}\mathcal H^{q,\mathrm{\acute{e}t}}_X(\mu_n^{\otimes q})$ 
for any smooth $k$-variety $X$ and an integer $n>0$, where $\mathcal K^M_{q,X}$ is restriction of 
$\mathcal K^M_q$ to $X$, $\mu_n$ is the \'etale sheaf of $n$-th roots of unity, and 
$\mathcal H^{q,\mathrm{\acute{e}t}}_X(\mu_n^{\otimes q})$ is the sheafification of the presheaf 
$U\mapsto H^q_{\mathrm{\acute{e}t}}(U,\mu_n^{\otimes q})$ for open $U\subseteq X$. 

Theorem~\ref{moving_lemma} implies easily (Proposition~\ref{sheafif}) that the \'etale and h- 
sheafifications of $\Omega^{q,\log}_{-|k}$ are $\Omega^{q,\log}_{-|k}\otimes\mathbb Q$ if $q\ge 1$. 

\section{Proofs} 
Recall (\cite{Voevodsky}) that Voevodsky's {\sl $(-1)$-construction} (or {\sl contraction}) 
associates to a Zariski presheaf $\mathcal F$ the presheaf $\mathcal F_{-1}:U\mapsto\mathcal F_{-1}(U):=
\mathrm{coker}\left[\mathcal F(U)\xrightarrow{p^*}\mathcal F(U\times\Gm)\right]$ for any 
smooth $U$,\footnote{or $\mathcal F_{-1}(U):=\mathrm{coker}\left[\mathcal F(\mathbb A^1_U)
\xrightarrow{j^*}\mathcal F(U\times\Gm)\right]$, where $U\times\Gm\xrightarrow{j}
\mathbb A^1_U:=U\times\mathbb A^1$ is the natural inclusion.} 
where $\Gm:=\mathop{\mathrm{Spec}}\nolimits\mathbb Z[T^{\pm 1}]$ is the multiplicative group 
scheme, and $U\times\Gm\xrightarrow{p}U$ is the (split) projection. 
Then $\mathcal F_{-n-1}:=(\mathcal F_{-n})_{-1}$ for all $n\ge 1$. A presheaf $\mathcal F$ 
is called {\sl homotopy invariant} if, for all $U$, the projections $\mathbb A^1_U\to U$ 
induce isomorphisms $\mathcal F(U)\xrightarrow{\sim}\mathcal F(\mathbb A^1_U)$. 

Another reminder: a {\sl presheaf with transfers} is an additive contravariant 
functor from an appropriate category of finite correspondences to the category of abelian groups. 
\begin{theorem}[\cite{Voevodsky}, Theorem 4.37; \cite{MVW}, Theorem 24.11] \label{main} Let $k$ be a perfect 
field, $X$ be a smooth scheme of dimension $n$ over $k$, and $\mathcal F$ a homotopy invariant Zariski sheaf 
with transfers. Then there is a canonical exact sequence of Zariski sheaves on $X$ of the form: 
\[0\to\mathcal F\to\coprod_{x\in X^{(0)}}i_{x*}\mathcal F\to\coprod_{x\in X^{(1)}}i_{x*}\mathcal F_{-1}
\to\cdots\to\coprod_{x\in X^{(n)}}i_{x*}\mathcal F_{-n}\to 0.\] 
\end{theorem} 

\subsection{The case of $\Omega^{\bullet,\log}_{-|k}$} To prove Theorem~\ref{log} by applying 
Voevodsky's Theorem~\ref{main} to $\Omega^{\bullet,\log}_{-|k}$, one has to \begin{itemize} 
\item construct a structure of a Zariski presheaf on $\Omega^{q,\log}_{-|k}$, 
\item check that $\Omega^{q,\log}_{-|k}$ with the constructed presheaf 
structure is in fact a Zariski sheaf, 
\item construct transfers, \item check homotopy invariance, and \item construct an 
isomorphism $\left(\Omega^{q,\log}_{-|k}\right)_{-1}=\Omega^{q-1,\log}_{-|k}$. \end{itemize} 

\begin{lemma} \label{Zariski-sheaf} $\Omega^{q,\log}_{-|k}:U\mapsto\Omega^{q,\log}_{-|k}(U)$ is a Zariski 
$($and even Nisnevich$)$ sheaf on the category of smooth irreducible $k$-varieties. \end{lemma} 
\begin{proof} If $\varphi:V\to U$ is a morphism, and 
$\eta=\sum_i\mathrm{d}\log f_{i1}\wedge\cdots\wedge\mathrm{d}\log f_{iq}\in\Omega^q_{U|k}(U)$ 
then $\varphi^*\eta\in\Omega^q_{V|k}(V)$ is well-defined. It remains to show that $\varphi^*\eta$ 
is a sum of exterior products of logarithmic differentials of rational functions on $V$. 

We choose a flag of subvarieties $V_s:=\mathrm{image}(\varphi)\subset\cdots\subset V_1\subset V_0\subseteq U$, 
where $V_i$ is of codimension $i$ in $U$, and proceed by induction on $s\ge 0$. 
As the blow-ups with smooth centres change neither function fields nor spaces of regular differentials, 
we may use Hironaka's resolution of singularities to assume that the $V_i$'s are smooth. 
The case $s=0$ (i.e. of dominant $\varphi$) is evident: 
$\varphi^*\eta=\sum_i\mathrm{d}\log(f_{i1}\circ\varphi)\wedge\cdots\wedge\mathrm{d}\log(f_{iq}\circ\varphi)$. 

For induction step, we may assume that the image of $\varphi$ is dense in an 
irreducuble divisor $D$ on $U$. Let $v$ be the discrete valuation of $k(U)$ 
corresponding to $D$, $\mO_v$ be the valuation ring of $v$, and $\pi$ be a uniformizer. 
Set $h_{ij}:=f_{ij}/\pi^{v(f_{ij})}$. Then 
\begin{multline}\label{loc-par}\eta=\sum_i\frac{\mathrm{d}h_{i1}}{h_{i1}}\wedge\cdots\wedge
\frac{\mathrm{d}h_{iq}}{h_{iq}}+\frac{\mathrm{d}\pi}{\pi}\wedge\omega,\ \mbox{ where }\\
\omega:=\sum_i\sum_{j=1}^q(-1)^{j-1}v(f_{ij})\frac{\mathrm{d}h_{i1}}{h_{i1}}\wedge\cdots\wedge
\widehat{\frac{\mathrm{d}h_{ij}}{h_{ij}}}\wedge\cdots\wedge\frac{\mathrm{d}h_{iq}}{h_{iq}}.\end{multline} 

Since $\omega\in\Omega^{q-1}_{\mO_v|k}$ and $\eta\in\Omega^q_{U|k}(U)\subset\Omega^q_{\mO_v|k}$, 
one has $\frac{\mathrm{d}\pi}{\pi}\wedge\omega\in\Omega^q_{\mO_v|k}$, and therefore, 
$\omega\in\pi\Omega^{q-1}_{\mO_v|k}+\mathrm{d}\pi\wedge\Omega^{q-2}_{\mO_v|k}$, 
so $\mathrm{d}\log\pi\wedge\omega$ should vanish on $D$, which means that 
the restriction of $\eta$ to $D$ is of desired type: 
$\eta|_D=\sum_i\frac{\mathrm{d}h_{i1}}{h_{i1}}\big\vert_D\wedge\cdots\wedge
\frac{\mathrm{d}h_{iq}}{h_{iq}}\big\vert_D$. \end{proof}

\begin{lemma}[Transfers] The transfer on Milnor $K$-theory $($constructed in \cite[Ch.I, \S5]{BT}$)$ 
descends to $\Omega^{q,\log}_{-|k}$, where it coincides with restriction of the transfer on differential 
$q$-forms $\Omega^q_{-|k}:V\mapsto\Omega^q_{V|k}(V)\ ($shown to preserve regularity in 
\cite{LecomteWach}$)$. In the case of transposed finite morphisms $U\to V$, it induces 
restriction to $\Omega^{q,\log}_{-|k}(U)$ of the trace $\Omega^q_{-|k}(U)\to\Omega^q_{-|k}(V)$. \end{lemma} 
\begin{proof} Let a finite correspondence from $U$ to $V$ be presented by an irreducible 
variety $W$, so $K:=k(W)$ is a finite field extension of $F:=k(U)$, and $k(V)\subseteq K$. Let 
$L|K$ be a finite field extensions such that $L|F$ is Galois. 
Let $K^M_q(F)\xrightarrow{j_{K|F*}}K^M_q(K)\xrightarrow{j_{L|K*}}K^M_q(L)$ be the maps induced by the 
inclusions of fields. Then the compositions $j_{K|F*}N_{L|F}j_{L|K*}$ and $[L:K]\cdot j_{K|F*}N_{K|F}$ 
coincide. But $j_{L|F*}N_{L|F}=\sum_{g\in\mathrm{Gal}(L|F)}g_*$ on $K^M_*(L)$, 
so the diagram \[\begin{diagram}\node{K^M_q(L)}\arrow[4]{e,t}{j_{L|F*}N_{L|F}=\sum_{g\in\mathrm{Gal}(L|F)}g_*}
\arrow{s,r}{\mathrm{d}\log}\node{}\node{}\node{}\node{K^M_q(L)}\arrow{s,r}{\mathrm{d}\log}\\
\node{\Omega^q_{L|k}}\arrow{e,t}{\cong}\node{L\otimes_F\Omega^q_{F|k}}
\arrow[2]{e,t}{\mathrm{tr}_{L|F}\cdot id=\sum_{g\in\mathrm{Gal}(L|F)}g_*\cdot id}
\node{}\node{\Omega^q_{F|k}}\arrow{e,t}{\subseteq}\node{\Omega^q_{L|k}}\end{diagram}\]
commutes, and therefore, $\mathrm{d}\log$ transforms $j_{K|F*}N_{K|F}$ to restriction of the 
composition 
$\Omega^q_{K|k}=K\otimes_F\Omega^q_{F|k}\xrightarrow{\mathrm{tr}_{K|F}\cdot id}\Omega^q_{F|k}
\hookrightarrow\Omega^q_{K|k}$. \end{proof} 

\begin{lemma} \label{H0-log-forms} For any smooth proper $k$-variety $X$ and an integer $q>0$, 
one has $\Omega^{q,\log}_{-|k}(X)=0$. \end{lemma} 
\begin{proof} Denote by $F^{\bullet}$ the Hodge filtration on the de Rham cohomology. 
The inclusions of complexes $\Omega^{q,\log}_{-|k}[-q]
\hookrightarrow\Omega^{\ge q}_{-|k}\hookrightarrow\Omega^{\bullet}_{-|k}$ induce inclusions 
$\Omega^{q,\log}_{-|k}(X)\subseteq F^qH^q_{\mathrm{dR}|k}(X)\subseteq H^q_{\mathrm{dR}|k}(X)$. 
Let $\eta=\sum_i\frac{\mathrm{d}f_{i1}}{f_{i1}}\wedge\cdots\wedge\frac{\mathrm{d}f_{iq}}{f_{iq}}
\in\Omega^{q,\log}_{-|k}(X)$, $U$ be the complement in $X$ to the union of the supports of 
the divisors of $f_{ij}$'s, and $\eta'\in F^qH^q_{\mathrm{dR}|k}(U)$ be the image of $\eta$ 
under the restriction $H^q_{\mathrm{dR}|k}(X)\xrightarrow{\rho}H^q_{\mathrm{dR}|k}(U)$. 

We apply the Lefschetz principle: together with $X$, all $f_{ij}$'s are defined over a finitely 
generated subfield of $k$, so we may assume that $k$ is finitely generated, and therefore, admits 
an embedding into $\mathbb C$. Replacing $X$ by $X\times_k\mathbb C$, we may further assume that 
$k=\mathbb C$, thus, identifying the algebraic de Rham cohomologies $H^q_{\mathrm{dR}|k}$ of $X$ 
and $U$ with the smooth de Rham and singular cohomologies, both with $\mathbb C$-coefficients, 
of the smooth manifolds $X(\mathbb C)$ and $U(\mathbb C)$. 

As $\frac{\mathrm{d}\overline{f}}{\overline{f}}=\mathrm{d}\log|f|^2-\frac{\mathrm{d}f}{f}$, 
the difference $\eta'-(-1)^q\overline{\eta'}$ is an exact $q$-form on $U(\mathbb C)$. 
This means that $\eta-(-1)^q\overline{\eta}$ belongs to the kernel of $\rho$, in particular, to 
the complexification of the coniveau-1 Hodge substructure of $H:=H^q(X(\mathbb C),\mathbb Z)$. 
However, $F^q(N^1H\otimes\mathbb C)=0$, and therefore, 
$\eta=(-1)^q\overline{\eta}\in F^q(H\otimes\mathbb C)\cap\overline{F^q}(H\otimes\mathbb C)$. 
But $F^q(H\otimes\mathbb C)\cap\overline{F^q}(H\otimes\mathbb C)=0$, 
since $H$ is a pure Hodge structure of weight $q$. \end{proof}

\begin{proof}[Proof of Corollary~\ref{free}] Proceed by induction on $q\ge 0$, the case $q=0$ being 
trivial. According to Nagata and Hironaka, any smooth connected $k$-variety $U_{\alpha}$ can be 
presented as an open subset of a smooth proper connected $k$-variety $X$. 
By Theorem~\ref{log}, there is an exact sequence 
$\Omega^{q,\log}_{-|k}(X)\to\Omega^{q,\log}_{k(X)|k}\to\bigoplus_{x\in X^{(1)}}\Omega^{q-1,\log}_{-|k}(x)$. 
By Lemma~\ref{H0-log-forms}, $\Omega^{q,\log}_{-|k}(X)=0$, so 
$\Omega^{q,\log}_{k(X)|k}$ is a subgroup of $\bigoplus_{x\in X^{(1)}}\Omega^{q-1,\log}_{-|k}(x)$, 
which is free by induction hypothesis, and therefore, $\Omega^{q,\log}_{k(X)|k}$ is free. 
Thus, $\mathop{\underrightarrow{\lim}}\limits_{\alpha}\Omega^{q,\log}_{-|k}(U_{\alpha})
\subseteq\Omega^{q,\log}_{k(X)|k}$ is free as well. \end{proof} 

\begin{proof}[Proof of Corollary~\ref{log-mod-n}] Let $A^{\bullet}$ be the Gersten 
resolution of $\Omega^{q,\log}_{-|k}$. By Corollary~\ref{free}, the sequence of complexes 
$0\to A^{\bullet}\xrightarrow{n\cdot}A^{\bullet}\to A^{\bullet}/n\to 0$ is exact. As 
$A^{\bullet}$ is acyclic in positive degrees, the corresponding  long exact cohomological 
sequence is $0\to\Omega^{q,\log}_{-|k}\xrightarrow{n\cdot}\Omega^{q,\log}_{-|k}\to ?\to 0$, 
so the complex $A^{\bullet}/n$ is also acyclic in positive degrees and 
$H^0(A^{\bullet}/n)=\Omega^{q,\log}_{-|k}/n$. \end{proof} 

\subsection{The case of $\Omega^{\bullet,\log}_{M,-|k}$} 
For any local factorial domain $A\supseteq k$ with the fraction field $F$, let $T(A)$ be 
the subgroup of $F(t)^{\times}$ generated by $A^{\times}$ and the monic polynomials in $A[t]$. 

Let $K_n^t(A)$ be the quotient of the $n$-th tensor power $T(A)^{\otimes n}$ of $T(A)$ by the subgroup 
generated by the tensors $p_1\otimes\cdots\otimes p_n$, where either $p_i+p_{i+1}=1$ or $p_{i+1}=-p_i$ 
for some $i$. The image of $p_1\otimes\cdots\otimes p_n$ in $K_n^t(A)$ is denoted $\{p_1,\dots,p_n\}$. 
Define the homomorphism $K_n^t(A)\xrightarrow{\partial_{\infty}}K_n^M(A)/\{-1\}\cdot K_{n-1}^M(A)$ 
on the generators by $\partial_{\infty}\{p_1,\dots,p_n\}:=
\sum_{j=1}^n(-1)^j\deg p_j\cdot\{\gamma_1,\dots,\widehat{\gamma_j},\dots,\gamma_n\}$, 
where $\gamma_j$ denotes the leading coefficient of $p_j$. 

\begin{theorem}[\cite{Kerz}, Theorem 4.4] \label{Kerz4.4} The Bass--Milnor--Tate short exact sequence 
in Milnor $K$-theory of fields 
$0\to K_n^M(F)\to K_n^M(F(t))\to\bigoplus_{\pi'}K_{n-1}^M(F[t]/(\pi'))\to 0$ restricts to an 
exact sequence $0\to K_n^M(A)\to K_n^t(A)\xrightarrow{\partial}\bigoplus_{\pi}K_{n-1}^M(A[t]/(\pi))\to 0$, 
where the direct sums are over all monic, irreducible $\pi'\in F[t]$ and $\pi\in A[t]$, 
respectively. \end{theorem} 

\begin{corollary}[Transfers II] The transfer on $\Omega^{q,\log}_{-|k}$ restricts 
to the transfer on $\Omega^{q,\log}_{M,-|k}$. \end{corollary} 
\begin{proof} One has to show that for any finite extension $A\subset B$ of local factorial 
domains the trace $\mathrm{tr}_{B|A}:\Omega^q_{B|k}\to\Omega^q_{A|k}$ restricts to a map 
$\Omega^{q,\log}_{M,B|k}\to\Omega^{q,\log}_{M,A|k}$, i.e. for any $b_j\in B^{\times}$ 
there exist $a_{ij}\in A^{\times}$ such that $\mathrm{tr}_{B|A}\eta=
\sum_i\frac{\mathrm{d}a_{i1}}{a_{i1}}\wedge\cdots\wedge\frac{\mathrm{d}a_{iq}}{a_{iq}}$, 
where $\eta=\frac{\mathrm{d}b_1}{b_1}\wedge\cdots\wedge\frac{\mathrm{d}b_q}{b_q}$. 
There is a monic polynomial $\pi\in A[t]$ such that $B=A[t]/(\pi)$. By Theorem~\ref{Kerz4.4}, 
there is an element $\alpha\in K_{q+1}^t(A)$ such that $\partial\alpha=\{b_1,\dots,b_q\}$. 
Then $\mathrm{tr}_{B|A}\eta=-\sum_i\mathrm{d}\log\partial_{\infty}
\alpha\in\Omega^{q,\log}_{M,A|k}$. \end{proof}

\begin{proof}[Proof of Theorem~\ref{moving_lemma}] The sheaf $\Omega^{q,\log}_{-|k}$ and 
its subsheaf $\Omega^{q,\log}_{M,-|k}$ have the same Gersten resolution, so they coincide. \end{proof} 

\begin{remark} Lemma~\ref{H0-log-forms} has the following generalization. {\it For any 
smooth proper $k$-variety $X$ and integers $q>p\ge 0$, the natural inclusion of the one-term 
complex $\Omega^{q,\log}_{-|k}[-q]$ into the de Rham complex $\Omega^{\bullet}_{-|k}$ 
induces the zero map} $H^p(X,\Omega^{q,\log}_{-|k})\to H^{p+q}_{\mathrm{dR}|k}(X)$. 
\begin{proof} We apply the Lefschetz principle again: together with $X$, any \v{C}ech cocycle 
presenting a class in $H^p(X,\Omega^{q,\log}_{M,-|k})$ is defined over a finitely generated 
subfield of $k$, so we may assume that $k$ is finitely generated, and therefore, admits 
an embedding into $\mathbb C$. The injectivity of the scalar extension from $k$ to 
$\mathbb C$ on de Rham cohomology allows us to replace $k$ by $\mathbb C$. 

For each integer $b\ge 0$, let $\mathcal H^{b,\mathrm{sing}}_{-(\mathbb C),\mathbb Z(q)}$ 
be the Zariski sheaf, associated to the presheaf of singular cohomology 
$U\mapsto H^b(U(\mathbb C),\mathbb Z(q))$. By \cite{Hodge-II}, it takes values in 
the category of (colimits of) Hodge structures with weights in the segment 
$[b-2q,2(b-q)]$. Therefore, the terms of the Bloch--Ogus (coniveau) 
spectral sequence $E^{a,b}_2=H^a(X,\mathcal H^{b,\mathrm{sing}}_{-(\mathbb C),\mathbb Z(q)})
\Rightarrow H^{a+b}(X(\mathbb C),\mathbb Z(q))$ carry Hodge structure, 
while the differentials are morphisms of Hodge structures.

The sheaf $\Omega^{q,\log}_{M,X|k}$ is naturally embedded into the sheaf 
$\mathcal H^{q,\mathrm{sing}}_{-(\mathbb C),\mathbb Z(q)}$ as a sheaf of pure 
Hodge structures of weight $0$. Then $H^p(X,\Omega^{q,\log}_{M,X|k})$ is also 
a pure Hodge structure of weight $0$. However, $H^{p+q}(X(\mathbb C),\mathbb Z(q))$ 
is a pure Hodge structure of weight $p-q<0$. \end{proof} \end{remark}

\begin{lemma}[Homotopy invariance and contraction] 
Let $\mathcal F_q$ be either $\Omega^{q,\log}_{-|k}$ or $\Omega^{q,\log}_{M,-|k}$. For any smooth 
$k$-variety $U$, the pullback $\mathcal F_q(U)\to\mathcal F_q(\mathbb A^1_U)$ is bijective. 
One has $\left(\mathcal F_q\right)_{-1}=\mathcal F_{q-1}$. \end{lemma} 
\begin{proof} One has $\Omega^q_{-|k}(\mathbb G_{\mathrm{m},U})=\Omega^q_{-|k}(U)[T^{\pm 1}]\oplus
\mathrm{d}T\wedge\Omega^{q-1}_{-|k}(U)[T^{\pm 1}]$, where $\mathbb G_{\mathrm{m},U}:=U\times\Gm$, 
and $T$ is a global invertible coordinate on $\Gm$, so 
\begin{equation}\label{sum-q-form}\Omega^q_{-|k}(\mathbb G_{\mathrm{m},U})=
\Omega^q_{-|k}(U)[T^{\pm 1}]+\mathrm{d}\log T\wedge\Omega^{q-1}_{-|k}(U)+
\mathrm{d}\left(\Omega^{q-1}_{-|k}(U)[T^{\pm 1}]\right).\end{equation} 

Denote by $\Omega^{q,\mathrm{cl}}_{-|k}\subset\Omega^q_{-|k}$ the subsheaf of closed $q$-forms. 
Passing to closed forms in (\ref{sum-q-form}) gives 
$\Omega^{q,\mathrm{cl}}_{-|k}(\mathbb G_{\mathrm{m},U})=
\Omega^{q,\mathrm{cl}}_{-|k}(U)\oplus\mathrm{d}\log T\wedge\Omega^{q-1}_{\mathrm{cl},-|k}(U)
\oplus\mathrm{d}((T-a)\Omega^{q-1}_{-|k}(U)[T^{\pm 1}])$ for any $a\in k^{\times}$ by 
sending $\eta$ to the triple $(p_a^*\eta,P_q\eta,\eta-p_a^*\eta-P_q\eta)$. Here $p_a^*$ and 
$P_q$ are orthogonal projectors on $\Omega^q_{-|k}(\mathbb G_{\mathrm{m},U})$ preserving 
$\Omega^{q,\mathrm{cl}}_{-|k}(\mathbb G_{\mathrm{m},U})$, where $p_a$ is the composition 
$\mathbb G_{\mathrm{m},U}\xrightarrow{\mathrm{pr}_U}U\xrightarrow{\times\{a\}}
\mathbb G_{\mathrm{m},U}$ and $P_q:=\frac{\mathrm{d}T}{T}\wedge\iota_{\xi}$ for 
the contraction $\iota_{\xi}$ with the vector field $\xi:=T\frac{\partial}{\partial T}$. 

If a closed form on $\mathbb G_{\mathrm{m},U}$ is not in $\Omega^{q,\mathrm{cl}}_{-|k}(U)+
\mathrm{d}\log T\wedge\Omega^{q-1,\mathrm{cl}}_{-|k}(U)$ then it has a pole of order $\ge 2$ 
along at least one of the divisors $T=0$ and $T=\infty$, so it cannot be logarithmic. 

As it shows equation (\ref{loc-par}) in the proof of Lemma~\ref{Zariski-sheaf} with $\pi=T$, if 
$\eta\in\mathcal F_q(\mathbb G_{\mathrm{m},U})$ then\footnote{If $\eta=\sum_i
\mathrm{d}\log f_{i1}^{(x)}\wedge\cdots\wedge\mathrm{d}\log f_{iq}^{(x)}$, 
where $f_{ij}^{(x)}\in\mathcal O^{\times}_{\mathbb G_{\mathrm{m},U},x}$, then there is 
an element $a\in\mathbb Q\subseteq k$ such that the divisor $T=a$ is contained the support of 
the divisor of neither of $f_{ij}^{(x)}$'s, so $f_{ij}^{(x)}|_{U\times\{a\}}$ is well-defined.} 
$\eta|_{U\times\{a\}}\in\mathcal F_q(U)$ is 
independent of $a$, and $\iota_{\xi}\eta\in\mathcal F_{q-1}(U)$. 

Thus, the logarithmic $q$-forms on $\mathbb G_{\mathrm{m},U}$ are contained in 
$\Omega^{q,\mathrm{cl}}_{-|k}(U)+\mathrm{d}\log T\wedge\Omega^{q-1,\mathrm{cl}}_{-|k}(U)$, and therefore, 
$\mathcal F_q(\mathbb G_{\mathrm{m},U})=\mathcal F_q(U)\oplus\mathrm{d}\log T\wedge\mathcal F_{q-1}(U)$. 
In particular, since \[\Omega^q_{-|k}(\mathbb A^1_U)=\Omega^q_{-|k}(U)[T]\oplus
\mathrm{d}T\wedge\Omega^{q-1}_{-|k}(U)[T]=\Omega^q_{-|k}(U)[T]+\mathrm{d}(\Omega^{q-1}_{-|k}(U)[T]),\] 
the logarithmic $q$-forms on $\mathbb A^1_U$ are contained in $\Omega^q_{\mathrm{cl},-|k}(U)$. \end{proof} 

\begin{proof}[Proof of Corollary~\ref{alg-closed-mod-n}] The composition 
$\Omega^{q,\log}_{X|k}\hookrightarrow\mathcal H^q_{\mathrm{dR},X|k}\hookrightarrow
\mathcal H^q_{\mathrm{dR},X(\mathbb C)}\cong\mathcal H^{q,\mathrm{sing}}_{X(\mathbb C),\mathbb C}$ 
factors through $\mathcal H^{q,\mathrm{sing}}_{X(\mathbb C),\mathbb Z(q)}$. Let $\Omega^{q,\log}_{X|k}
\xrightarrow{\beta}\mathcal H^{q,\mathrm{sing}}_{X(\mathbb C),\mathbb Z(q)}$ be this composition. 
As above, denote by $\mu_n$ the sheaf of $n$-th roots of unity. 
Obviously, there is a unique morphism $\Omega^{q,\log}_{X|k}/n\xrightarrow{\gamma}
\mathcal H^{q,\mathrm{\acute{e}t}}_{X\times_k\mathbb C}(\mu_n^{\otimes q})$ making 
commutative the diagram\\ \scalebox{0.97}{$
\begin{diagram} \node{\mathcal K^M_{q,X}}\arrow{e,t}{\mathrm{d}\log}\arrow{s,r}{\bmod n}
\node{\Omega^{q,\log}_{X|k}}\arrow{e,t}{\beta}\arrow{s,r}{\bmod n} 
\node{\mathcal H^{q,\mathrm{sing}}_{X(\mathbb C),\mathbb Z(q)}}\arrow{e,t}{\bmod n} 
\node{\mathcal H^{q,\mathrm{sing}}_{X(\mathbb C),\mathbb Z(q)}/n}\arrow{e,t}{\hookrightarrow} 
\node{\mathcal H^{q,\mathrm{sing}}_{X(\mathbb C),\mathbb Z/n(q)}}\arrow{s,r}{\cong}\\ 
\node{\mathcal K^M_{q,X}/n}\arrow{e,t}{\mathrm{d}\log_n}\node{\Omega^{q,\log}_{X|k}/n}\arrow[3]{e,t}{\gamma} 
\node{}\node{}\node{\mathcal H^{q,\mathrm{\acute{e}t}}_{X\times_k\mathbb C}(\mu_n^{\otimes q})}\end{diagram}$} 
Clearly, $\mathrm{d}\log_n$ is surjective. The lower horizontal composition coincides with the 
composition $\mathcal K^M_{q,X}\to\mathcal H^{q,\mathrm{\acute{e}t}}_X(\mu_n^{\otimes q})\xrightarrow{\sim}
\mathcal H^{q,\mathrm{\acute{e}t}}_{X\times_k\mathbb C}(\mu_n^{\otimes q})$ (since $k$ is algebraically 
closed), so it is injective by the Bloch--Kato conjecture proved in 
\cite{V_Bloch-Kato}. This means that $\mathrm{d}\log_n$ is an isomorphism. \end{proof} 

\section{A remark on sheafifications} 
\begin{definition}[\cite{Voevodsky96}, Definitions 3.1.1, 3.1.2] 
A morphism of schemes $X\xrightarrow{p}Y$ is called a {\sl universal topological epimorphism} if 
its base change $Z\times_YX\xrightarrow{p_f}Z$ by any morphism of schemes $Z\xrightarrow{f}Y$ 
is a topological epimorphism, i.e. if (i) $p_f$ is surjective and (ii) a subset $A$ in $Z$ is 
Zariski open if and only if the subset $p_f^{-1}(A)$ is Zariski open in $Z\times_YX$. 

The h-{\sl topology} on the category of separated $k$-schemes of finite type is the Grothendieck 
topology associated with the pretopology whose coverings are finite families of the form 
$\{U_i\xrightarrow{p_i}X\}_{i\in I}$ such that the morphism $\coprod_{i\in I}p_i:\coprod_{i\in I}U_i\to X$ 
is a universal topological epimorphism. \end{definition}

It is shown in \begin{itemize} \item \cite[Proposition 4.2]{Lee} and \cite[Theorem 3.6]{HuberJoerder} 
that $\Omega^q_{-|k}$ is an h-sheaf on the category of smooth $k$-varieties; 
\item \cite[Proposition 2.2]{Scho12} that any h-sheaf has a canonical transfer structure. \end{itemize} 

\begin{proposition} \label{sheafif} If $q\ge 1$ then $\Omega^{q,\log}_{-|k}\otimes\mathbb Q$ is both 
\'etale and h- sheafification of $\Omega^{q,\log}_{-|k}$. \end{proposition} 
\begin{proof} For each integer $n>1$, set 
$G_n:=\mathop{\mathrm{Spec}}\nolimits(k[t_1^{\pm 1/n},\dots,t_q^{\pm 1/n}])$, so 
the inclusion of the coordinate rings induces an  \'etale covering $G_n\to G:=G_1$. 
Then $\Omega^{q,\log}_{-|k}(G_n)=\frac{1}{n^q}\mathbb Z\cdot\frac{\mathrm{d}t_1}{t_1}
\wedge\cdots\wedge\frac{\mathrm{d}t_q}{t_q}$, so $\Omega^{q,\log}_{-|k}(G_n)\xrightarrow{id}
\Omega^{q,\log}_{-|k}(G_n)\rightrightarrows\Omega^{q,\log}_{-|k}(G_n\times_GG_n)$ is the equalizer 
diagram.\footnote{If $u_1^n=u_2^n\in\mO(G)^{\times}\subset\mO(G_n\times_GG_n)^{\times}$ 
for some $u_i\in\mathrm{pr}_i^*\mO(G_n)^{\times}\subset\mO(G_n\times_GG_n)^{\times}$ 
($i=1,\ 2$) then $du_1/u_1=du_2/u_2\in\Omega^{1,\log}_{-|k}(G_n\times_GG_n)^{\times}$.} 

Denote by $\Omega^{q,\log}_{\mathrm{\acute{e}t},-|k}$ the \'etale sheafification of $\Omega^{q,\log}_{-|k}$. 
Then $\Omega^{q,\log}_{-|k}(G)\otimes\mathbb Q\subseteq\Omega^{q,\log}_{\mathrm{\acute{e}t},-|k}(G)$. 

Let $X$ be a smooth irreducible $k$-scheme and 
$\eta=\sum_i\frac{\mathrm{d}f_{i1}}{f_{i1}}\wedge\cdots\wedge\frac{\mathrm{d}f_{iq}}{f_{iq}}
\in\Omega^{q,\log}_{-|k}(X)$. For any $x\in X$, denote by $\mO_x\subset k(X)$ the local ring of $x$. 

Since $\Omega^{q,\log}_{-|k}=\Omega^{q,\log}_{M,-|k}$, for each $x\in X$, there exists 
a collection $(f^{(x)}_{ij}\in\mO_x^{\times})$ such that $\eta=\sum_i
\frac{\mathrm{d}f^{(s)}_{i1}}{f^{(s)}_{i1}}\wedge\cdots\wedge\frac{\mathrm{d}f^{(s)}_{iq}}{f^{(s)}_{iq}}$.  
Let $U_x\ni x$ be a Zariski open subset of $X$ such that $f^{(x)}_{ij}\in\mO^{\times}(U_x)$ 
for all $i,j$. 
As, for each $x$, the $q$-tuples $f^{(x)}_i:=\left(f^{(x)}_{i1},\dots,f^{(x)}_{iq}\right)$ can be 
considered as morphisms $U_x\to G$, this shows that $\Omega^{q,\log}_{\mathrm{\acute{e}t},-|k}(U_s)$ 
contains $f^{(x)*}_i\Omega^{q,\log}_{-|k}(G)\otimes\mathbb Q$. 
In particular, $\frac{1}{n^q}\eta\in\Omega^{q,\log}_{-|k}(X)$, and therefore, 
$\Omega^{q,\log}_{-|k}\otimes\mathbb Q\subseteq\Omega^{q,\log}_{\mathrm{\acute{e}t},-|k}$. 

As $\Omega^{q,\log}_{-|k}$ is a Zariski sheaf with transfers, the h-descent for 
$\Omega^{q,\log}_{-|k}\otimes\mathbb Q$ results from the following theorem of Voevodsky 
(\cite[Theorem 3.1.9]{Voevodsky96}): 

{\it Any h-covering $\{U_i\to X\}$ of an excellent reduced noetherian scheme $X$ admits a refinement 
$\{U_j'\xrightarrow{s\circ f\circ\mathrm{in}_j}X\}$, where $\{\mathrm{in}_j:U_j'\to\bar{U}\}$ 
is an open covering, $\bar{U}\xrightarrow{f}X_Z$ is a finite surjective morphism and 
$X_Z\xrightarrow{s}X$ is the blowup of a closed subscheme in $X$.} \end{proof}

\begin{remark} One can modify the argument of \cite[Theorem 1]{HuberJoerder} 
to show that the h-shea\-fi\-fi\-ca\-tion of $\Omega^{q,\log}_{-|k}$ extends to a subsheaf 
$\Omega^{q,\log}_{\mathrm{h},-|k}(X)$ of the h-sheaf $\Omega^q_{\mathrm{h},-|k}(X)$ 
on the category of separated schemes of finite type over $k$: 
\[\Omega^{q,\log}_{\mathrm{h},-|k}(X)=\left\{(\alpha_f)_{f:Y\to X}\in\prod\limits_{f:Y\to X,\ Y\text{ smooth}}
\Omega^{q,\log}_{-|k}(Y)\otimes\mathbb Q\ \vert\ \mbox{\scalebox{0.7}{$
\begin{diagram}\node{Y'} \arrow{s,l}{\phi} \arrow{se,t}{f'}\\ 
\node{Y} \arrow{e,t}{f} \node{X}\end{diagram}$}}\ \Rightarrow\ \phi^*\alpha_f=\alpha_{f'}\right\}.\] 
Or more explicitly, cf. \cite[Remark 3.8]{HuberJoerder}, $\Omega^{q,\log}_{\mathrm{h},-|k}(X)
=\mathrm{Ker}(\Omega^{q,\log}_{-|k}(X')\xrightarrow{\psi_1^*-\psi_2^*}\Omega^{q,\log}_{-|k}(X''))
\otimes\mathbb Q$, where $X'\to X$ and $X''\xrightarrow{\psi}X'\times_XX'$ are smooth h-covers. \end{remark}

\noindent 
{\sl Acknowledgements.} {\small The above proof of the Gersten property for logarithmic forms is 
a realization of a plan suggested to me by Vladimir Voevodsky at the IAS in Fall 2009. 

Questions of Anna Votinova helped me to clarify several points of this note. 

I am grateful to the Institute for Advanced Studies for its hospitality and 
to the National Science Foundation for support under agreement No. DMS-0635607.

Support from the Basic Research Program of HSE University is gratefully acknowledged (HSE-BR-2025-061). 

}


\vspace{5mm}

\begin{thebibliography}{} 
\bibitem{BT} H.Bass, J.Tate, The Milnor ring of a global field, 
Algebraic K-theory II (Battelle Memorial Inst., Seattle, Washington, 1972), 
Lecture Notes in Math., vol. 342, Springer-Verlag, New York, 1973, 349--446. 
\bibitem{Hodge-II} P.Deligne, Th\'eorie de Hodge II, Publ. Math. de l'IH\'ES, 40 (1971), 5--57. 
\bibitem{HuberJoerder} A.Huber, C.Jörder, Differential forms in the h-topology, 
Algebraic Geometry 1 (4) (2014) 449--478. 
\bibitem{Kerz} M.Kerz, The Gersten conjecture for Milnor K-theory, Invent. Math. 175 (2009) 1--33. 
\bibitem{LecomteWach} F.Lecomte, N.Wach, Le complexe motivique de De Rham, 
manuscripta math. 129, (2009) 75--90. 
\bibitem{Lee} B.Lee, Local acyclic fibrations and the de Rham complex. 
Homology, Homotopy Appl. 11:1 (2009), 115--140. 
\bibitem{MVW} C.Mazza, V.Voevodsky, C.Weibel, Lecture Notes on Motivic Cohomology, 
(Clay Mathematics Monographs Volume 2), 2006. 
\bibitem{Scho12} J.Scholbach, Geometric motives and the h-topology, Math. Z. 272 (2012), 965--986. 
\bibitem{Voevodsky96} V.Voevodsky, Homology of schemes, Selecta Math. (N.S.) 2, no. 1 (1996), 111--153. 
\bibitem{Voevodsky} V.Voevodsky, Cohomological theory of presheaves with transfers, Ann. Math. Studies, 
2000, 87--137. 
\bibitem{V_Bloch-Kato} V.Voevodsky, On motivic cohomology with $\mathbf Z/l$-coefficients, 
Ann. of Math. 174 (1) (2011), 401--438. 
\end{thebibliography}
\end{document}